\documentclass[11pt]{article}
\usepackage{graphicx}

\usepackage{amsmath}
\usepackage{epsfig}  
\usepackage{verbatim}  
\usepackage{tabularx}
\usepackage{xcolor}
\usepackage{lineno}
\usepackage{amssymb,latexsym,amscd,amsmath,amsfonts,enumerate,supertabular}
\usepackage{subfigure}
\usepackage{longtable}
\newtheorem{Theorem}{Theorem}[section] 

\newtheorem{lemma}[Theorem]{Lemma} 
\newtheorem{Proposition}[Theorem]{Proposition}

\begin{document}

\title{{\bf Numerical solution of a nonlinear functional integro-differential equation   }}

%
%
%
%
%

\author{ Dang Quang A$^{\text a}$,  Pham Huy Dien $^{\text a}$, Dang Quang Long$^{\text b}$\\
$^{\text a}$ {\it\small Center for Informatics and Computing, VAST}\\
{\it\small 18 Hoang Quoc Viet, Cau Giay, Hanoi, Vietnam}\\
{\small Email: dangquanga@cic.vast.vn, phdien@cic.vast.vn}\\
$^{\text b}$ {\it\small Institute of Information Technology, VAST,}\\
{\it\small 18 Hoang Quoc Viet, Cau Giay, Hanoi, Vietnam}\\
{\small Email: dqlong88@gmail.com}}
\date{ }
%


\maketitle

%

\date{}          

\begin{abstract}
In this paper, we consider  a boundary value problem (BVP) for a fourth order nonlinear functional integro-differential equation. We establish the existence and uniqueness of  solution and construct a numerical method for solving it. We prove that the method is of second order accuracy and obtain an estimate for the total error. Some examples demonstrate the validity of the obtained theoretical results and the efficiency of the numerical method.

\end{abstract}
{\small
}

\section{Introduction}

In this paper we consider the following problem for functional integro-differential equation (FIDE)
\begin{equation}\label{eq3}
\begin{split}
u^{(4)}(x)& =f(x,u(x),u(\varphi (x)), \int_0^1 k_0(x,t) u(t) dt, \int_0^1 k_1(x,t) u(\varphi(t)) dt ),\\
u(0)&=0 , \; u(1)=0, \; u''(0)=0 , \; u''(1)=0,
\end{split}
\end{equation}
where the function $f(x,u,y,v,z), k_0(x,t), k_1(x,t)$ and $\varphi (t)$ are assumed to be continuous, $\varphi : [0,1] \rightarrow [0,1]$. 
This is a natural combination or hybridization of the following problems:\\
{\bf Problem 1:} The problem for  functional differential equation
\begin{equation*}\label{eq2}
\begin{split}
x^{(2p)}(t)=f(t,x(t), x(\varphi (t))), \quad t\in [a,b],\\
x^{(i)}(a)=a_i, \; x^{(i)}(b)=b_i,\quad i=\overline{0, p-1}
\end{split}
\end{equation*}
where $\varphi : [a,b] \rightarrow \mathbb{R}, \; a \le \varphi (t) \le b, \forall t\in [a,b]$ and $f(t,x,y)$ are continuous functions.\\
{\bf Problem 2:} The problem for  integro-differential equation
\begin{equation*}\label{eq1}
\begin{split}
y^{(4)}(x)& =f(x,y(x), \int_0^1 k(x,t) y(t) dt ),\; 0<x<1,   \\
y(0)&=0 , \; y(1)=0, \; y''(0)=0 , \; y''(1)=0,
\end{split}
\end{equation*}
where the functions $f(x,y,v), k(x,t)$ are assumed to be continuous.\par
Problem 1 was considered in the work of Bica et al. \cite{Bica}. For solving this problem, the authors constructed successive approximations for the  equivalent integral equation with the use of cubic spline interpolation at each iterative step. The error estimate was obtained under strong conditions including the satisfaction of Lipschitz conditions of the function $f(x,y,v)$ in the unbounded domain $[a, b] \times \mathbb{R} \times \mathbb{R}$. Although in \cite{Bica} the method is constructed for the general function $\varphi (t)$ but in all numerical examples only the particular case $\varphi (t) =\alpha  t$ was considered and the conditions of convergence were not verified. It is a regret that in all examples  the Lipschitz condition for the function $f(x,y,v)$ is not satisfied. To overcome the limitations of Bica, very recently in \cite{Dang8} Dang and Dang proposed a new approach for the third order FDE. In this work the existence and uniqueness of solution and the error estimate of the approximate solution obtained by a discrete iterative method  are established under easily verified conditions. This approach is applicable to FDE of any order.\par
It should be said that functional differential equations (FDE)  have numerous applications in engineering and sciences \cite{Hale}. Therefore, for the last decades they have been studied by many authors. There are many works concerning the numerical solution of both initial and boundary value problems for them. The methods used are diverse including collocation method \cite{Reut,Yang}, iterative methods \cite{Bica,Khuri}, neural networks \cite{Hou,Raja}, and so on.\par
Problem 2 was considered in the last year, 2020, by Wang \cite{Wang}. Using the monotone method and a maximum principle, Wang constructed the sequences of functions, which converge to the extremal solutions of the problem. Very recently, Dang and Dang \cite{Dang9} studied a generalized version of the above problem, namely, the problem where in the right-hand side there is an extra $y'$. Using the method developed in the previous papers \cite{Dang1,Dang2,Dang4,Dang5,Dang6} they establish the existence and uniqueness of the solution and propose a solution method of  second order convergence. \par
As functional differential equations,  integro-differential equations (IDE) are the mathematical models of many phenomena in physics, biology, hydromechanics, chemistry, etc. In general, it is impossible to find the exact solutions of the problems involving these equations, especially when they are nonlinear. Therefore, many analytical approximate  methods and numerical methods have been developed for  these equations (see, e.g. \cite{Aruc,Chen,Chen1,Dang-Vu,Dasc,Fari,Lake,Sing,Swei,Tahe,Yula,Wang,Zhua}). \par
In the present paper, using the techniques of \cite{Dang8,Dang9} we investigate the existence and uniqueness of solutions of the problem \eqref{eq3} and construct iterative methods for solving the problem. We prove that the discrete iterative method is of second order accuracy. Some examples demonstrate the validity of the obtained theoretical results and the efficiency of the numerical method.

\section{Existence results}
To investigate the problem \eqref{eq3} for brevity we 
denote by $ \mathcal{K}_0, \mathcal{K}_1$ the integral operators defined by
$$(\mathcal{K}_0u)(x) = \int_0^1 k_0(x,t) u(t) dt, \; (\mathcal{K}_1u)(x) = \int_0^1 k_1(x,t) u(\varphi (t)) dt$$
and
set
\begin{equation}\label{eq4}
\psi (x)= f(x,u(x),u(\varphi (x)), (\mathcal{K}_0u)(x), (\mathcal{K}_1u)(x) ).
\end{equation}
Then the problem  becomes
\begin{equation}\label{eq5}
\begin{split}
u^{(4)}(x)& =\psi (x),\\
u(0)=0 , \; u(1)&=0, \; u''(0)=0 , \; u''(1)=0.
\end{split}
\end{equation}
It has a solution presented in the form
\begin{equation}\label{eq6}
u(x)=\int_0^1 G(x,s)\psi(s)ds, \quad 0<x<1,
\end{equation}
where $G(x,s)$ is the Green function of the operator $u''''(x)=0$ associated with the homogeneous boundary conditions $u(0)=u''(0)= u(1)=u''(1)=0$,
\begin{align}\label{eq7}
G(x,s)=\frac{1}{6}
\begin{cases}
s(x-1)(x^2-x+s^2), \quad & 0\leq s\leq x\leq 1 \\
x(s-1)(s^2-s+x^2), & 0\leq x\leq s\leq 1.
\end{cases}
\end{align}
Notice that $G(x,s) \ge 0 $ for all $(x,s) \in [0,1] \times [0,1]$.\\
It is obvious that the solution \eqref{eq6} depending on $\psi$ must satisfy the relations \eqref{eq4}. We shall express this by an operator equation for $\psi$. For this purpose, in the space $C[0, 1]$ introduce the operator $A$ by the formula
\begin{equation}\label{eq8}
(A\psi ) (x)= f(x,u(x),u(\varphi (x)), (\mathcal{K}_0u)(x), (\mathcal{K}_1u)(x) ),
\end{equation}
where $u(x)$ is the solution of the boundary value problem \eqref{eq5}.\\

It is easy to verify the following lemma.
\begin{lemma}
If the function $\psi$ is a fixed point of the operator $A$, i.e., $\psi$ is the solution of the operator equation 
\begin{equation}\label{eq9}
A\psi = \psi ,
\end{equation} 
where $A$ is defined by \eqref{eq8},
then the function $u(x)$ determined from the BVP \eqref{eq5} is a solution of the BVP \eqref{eq3}. Conversely, if the function $u(x)$ is the solution of the BVP \eqref{eq3} then the function $\psi (x)$ determined from \eqref{eq4}
satisfies the operator equation \eqref{eq9}.
\end{lemma}
Due to the above lemma we shall study the original BVP \eqref{eq3} via the operator equation \eqref{eq8}.

Set 
\begin{align}\label{eq10}
\begin{split}
M_0 & =\max_{0\leq x\leq 1} \int_0^1 |G(x,s)|ds ,\\
K_0 & =\max_{0\leq x\leq 1} \int_0^1 |k_0(x,s)|ds , \\
K_1 & =\max_{0\leq x\leq 1} \int_0^1 |k_1(x,s)|ds
\end{split}
\end{align}  
It is easy to obtain 
\begin{equation}\label{eq11}
 M_0=   \frac{5}{384} .
\end{equation}
Now for any positive number $M$, we define the domain
\begin{align}\label{eq12}
\begin{split}
\mathcal{D}_M = \{(x,u,y,v,z)\; | \; &0\leq x\leq 1,\; |u|\leq M_0M,\\
 &|y|\leq M_0M,\; |v| \le M_0K_0M,\;  |z|\leq M_0K_1M \}.
\end{split}
\end{align}
As usual, we denote by $B[0, M]$ the closed ball centered at $0$ with radius $M$ in the space $C[0,1]$, i.e.,
\begin{align*}
B[0, M] =\lbrace u\in C[0,1] \; \vert \; \|u\| \le M \rbrace ,
\end{align*}
where $\|u\| =\max _{0\le x \le 1} |u(x)|$.
\begin{Theorem}[Existence and uniqueness]\label{thm1}
Suppose that the functions $k_0(x,t)$, $k_1(x,t)$ are continuous in the square $[0,1] \times [0,1]$, the function $\varphi (t)$ is continuous and maps $[0, 1]$ to $[0,1]$, and there exist numbers $M>0$, $L_0,L_1,L_2,L_3\geq 0$ such that:
\begin{description}
\item [(i)]The function $f(x,u,y,v,z)$ is continuous in the domain $\mathcal{D}_M$ and\\
 $|f(x,u,y,v,z)|\leq M$, $ \forall (x,u,y,v,z)\in \mathcal{D}_M$.
\item [(ii)]$|f(x_2,u_2,y_2,v_2,z_2)-f(x_1,u_1,y_1, v_1,z_1)| \leq L_0 |u_2-u_1|+ L_1 |y_2-y_1|+ L_2 |v_2-v_1| +|L_3 |z_2-z_1|, \; \forall (x_i,u_i,y_i,v_i,z_i)\in \mathcal{D}_M, \; i=1,2.$
\item [(iii)]$q = (L_0+L_1+L_2K_0+L_3K_1)M_0 <1$.
\end{description}
Then the problem \eqref{eq3} has a unique solution $u\in C^4[0,1]$ satisfying $|u(x)| \leq M_0M$  for any $0 \le x \le 1.$
\end{Theorem}
\noindent {\bf Proof.}
Under the assumptions of the theorem we shall prove that the operator $A$ is a contraction mapping in the closed ball $B[O,M]$. Then the operator equation \eqref{eq9} has a unique solution $u \in C^{(4)}[0,1]$ and this implies the existence and uniqueness of solution of the BVP \eqref{eq3}.\par 
Indeed, take $\psi \in B[O,M]$. Then the problem \eqref{eq5} has a unique solution of the form \eqref{eq6}. From there and \eqref{eq10} we obtain $|u(x)| \le M_0\| \psi \|$ for all $x\in [0, 1]$. 
Further, from the last equation in \eqref{eq10} we have the estimate $| (\mathcal{K}_i u)(x)| \le M_0 K_i \|\psi\|$, $ x\in [0,1]$. Thus, if $\psi \in B[O,M]$, i.e., $\|\psi\| \le M$ then for any $x \in [0,1]$ we have
$$ |u(x)| \le M_0M, \;  |(\mathcal{K}_i u) (x)|\le M_0K_iM \; (i=0,1).
$$
Therefore, $(x, u(x),u(\varphi(x)), (\mathcal{K}_0 u )(x)), (\mathcal{K}_1 u)  (x))\in \mathcal{D}_M$. By the assumption $(i)$ there is
$$ |f(x, u(x),u(\varphi(x)), \mathcal{K}_0 u (x)), \mathcal{K}_1 u  (x)| \le M \quad \forall x\in [0,1].
$$
Hence, $|(A\psi)(x)| \le M, \; \forall x\in [0,1]$ and $\|A\psi \| \le M$. It means that $A$ maps $B[O,M]$ into itself.\par 
Next, take $\psi _1, \psi_2 \in B[O,M]$. Using the assumption $(ii)$ and $(iii)$ it is easy to obtain
\begin{align*}
\|A\psi _2 - A\psi _1\| \le ((L_0+L_1+L_2K_0+L_3K_1)M_0 )\|\psi _2 -\psi_1 \| =q \|\psi _2 -\psi_1 \|.
\end{align*}
Since $q<1$ the operator $A$ is a contraction map in $B[O,M]$. This completes the proof of the theorem.

Now, in order to study positive solutions of the BVP \eqref{eq3} we introduce the domain
\begin{align}\label{eq:A10}
\begin{split}
\mathcal{D}_M^+ = \{(x,u,y,v,z)\; | \; &0\leq x\leq 1,\; 0\le u\leq M_0M,\;  |y|\leq M_0M,\;\\
 & |v|\leq M_0K_1M,\; |z|\leq M_0K_2M \}.
\end{split}
\end{align}
and denote
$$ S_M = \{\psi \in C[0,1], 0 \le \psi (x) \le M   \}.
$$
\begin{Theorem}[Positivity of solution]\label{thm2}
If  in Theorem \ref{thm1} the condition $(i)$ is replaced by the condition\\
$(i'): \; 0\le f(x,u,y,v,z)\leq M, \; \forall (x,u,y,v,z)\in \mathcal{D}_M^+$ and $f(x,0,0,0,0) \not \equiv 0$,\\
and in the condition $(ii)$ replace $\mathcal{D}_M$ by $\mathcal{D}_M^+$,
then the problem \eqref{eq3} has a unique positive solution $u\in C^4[0,1]$ satisfying $0 \le u(x) \leq M_0M$  for any $0 \le x \le 1.$
\end{Theorem}
\noindent {\bf Proof.}
Similarly to the proof of Theorem \ref{thm1}, where instead of $\mathcal{D}_M$ and $B[O; M]$
there stand $\mathcal{D}_M^+$ and $S_M$, we conclude that the problem  has a nonnegative solution. Due to the condition $f(x,0,0,0,0) \not \equiv 0$, this solution must be positive.


\section{Numerical method }\label{IterMeth}

In this section we suppose that all the conditions of Theorem \ref{thm1} are satisfied. Then the problem \eqref{eq3} has a unique solution. For finding this solution consider the following iterative method:
\begin{enumerate}
\item Given 
\begin{equation}\label{iter1c}
\psi_0(x)=f(x,0,0,0,0).
\end{equation}
\item Knowing $\psi_m(x)$  $(m=0,1,...)$ compute
\begin{equation}\label{iter2c}
\begin{split}
u_m(x) &= \int_0^1   G(x,t)\psi_m(t)dt  ,\\
y_m(x) &:= u_m (\varphi (x))= \int_0^1 G(\varphi (x),t)\psi_m(t)dt ,\\
v_m(x) &= \int_0^1 k_0(x,t)u_m(t) dt,\\
z_m(x) &= \int_0^1 k_1(x,t) y_m(t) dt,\\
\end{split}
\end{equation}
\item Update
\begin{equation}\label{iter3c}
\begin{split}
\psi_{m+1}(x) &= f(x,u_m(x),y_m(x),v_m(x),z_m(x)).
\end{split}
\end{equation}
\end{enumerate}
This iterative method indeed is the successive iterative method for finding the fixed point of the operator $A$. Therefore, it converges with the rate of geometric progression and there holds the estimate
$$\|\psi_m -\psi\| \leq \frac{q^m}{1-q}\|\psi_1 -\psi_0\| = p_m d,$$
where  $\psi$ is the fixed point of the operator $A$ and
\begin{equation}\label{eq:pd}
p_m=\tfrac{q^m}{1-q},\; d=\|\psi_1 -\psi_0\|.
\end{equation}

This estimate implies the following result of the convergence of the iterative method \eqref{iter1c}-\eqref{iter3c}.
\begin{Theorem}\label{thm3}
Under the conditions of Theorem \ref{thm1} the iterative method \eqref{iter1c}-\eqref{iter3c} converges and for the approximate solution $u_m(t)$ there hold estimates
\begin{align*}
\|u_m-u\| &\leq M_0p_md, 
\end{align*}
where $u$ is the exact solution of the problem \eqref{eq3}, $p_m$ and $d$ are defined by \eqref{eq:pd}.
\end{Theorem}

To numerically realize the above iterative method we construct a corresponding discrete iterative method. For this purpose cover the interval $[0, 1]$   by the uniform grid $\bar{\omega}_h=\{x_i=ih, \; h=1/N, i=0,1,...,N  \}$ and denote by $\Psi_m(x), U_m(x),  Y_m(x), V_m(x), Z_m(x)$ the grid functions, which are defined on the grid $\bar{\omega}_h$ and approximate the functions $\psi (x), u_m(x),  y_m(x), v_m(x), z_m(x)$ on this grid. \par
Consider now the following discrete iterative method:
\begin{enumerate}
\item Given 
\begin{equation}\label{iter1d}
\Psi_0(x_i)=f(x_i,0,0,0,0),\ i=0,...,N. 
\end{equation}
\item Knowing $\Psi_m(x_i), \; m=0,1,...; \; i=0,...,N, $  compute approximately the definite integrals \eqref{iter2c} by the trapezium formulas
\begin{equation}\label{iter2d}
\begin{split}
U_m(x_i) &= \sum _{j=0}^N h\rho_j G(x_i,x_j)\Psi_m(x_j) ,\\
Y_m(x_i) &= \sum _{j=0}^N h\rho_j G( \varphi(x_i),x_j)\Psi_m(x_j) ,\\
V_m(x_i) &= \sum _{j=0}^N h\rho_j k_0(x_i,x_j)U_m(x_j) ,\\
Z_m(x_i) &= \sum _{j=0}^N h\rho_j k_1(x_i,x_j)Y_ m(x_j) ,\;  i=0,...,N,
\end{split}
\end{equation}
\noindent where $\rho_j$ is the weight of the trapezium formula, namely
\begin{equation*}
\rho_j = 
\begin{cases}
1/2,\; j=0,N\\
1, \; j=1,2,...,N-1.
\end{cases}
\end{equation*}

\item Update
\begin{equation}\label{iter3d}
\Psi_{m+1}(x_i) = f(x_i,U_m(x_i),Y_m(x_i), V_m(x_i),Z_m(x_i) ).
\end{equation}
\end{enumerate}
\noindent {\bf Remark 1.} 
In the case if the boundary conditions in the problem \eqref{eq3} are non-homogeneous, namely, if 
\begin{equation*}
u(0)=c_1 , \; u(1)=c_2, \; u''(0)=c_3 , \; u''(1)=c_4,
\end{equation*}
then the solution of the equation $u^{(4)}(x) =\psi (x)$ satisfying these boundary conditions has the form
\begin{equation}
u(x)=\int_0^1 G(x,s)\psi(s)ds +p(x), \quad 0<x<1.
\end{equation}
Here 
\begin{equation}
p(x)=a_0+a_1x+a_2x^2+a_3x^3
\end{equation}
with
\begin{align*}
a_0&=c_1, \; a_1=-c_1 +c_2-\frac{1}{3}c_3 -\frac{1}{6}c_4,\\
a_2&=\frac{1}{2}c_3, \; a_3= \frac{1}{6}c_4 -\frac{1}{6}c_3.
\end{align*}
Then the formulas in \eqref{iter2d} have the form
\begin{equation}\label{iter2da}
\begin{split}
U_m(x_i) &= \sum _{j=0}^N h\rho_j G(x_i,x_j)\Psi_m(x_j)+p(x_i) ,\\
Y_m(x_i) &= \sum _{j=0}^N h\rho_j G( \varphi(x_i),x_j)\Psi_m(x_j) +p(\varphi(x_i)) ,\\
V_m(x_i) &= \sum _{j=0}^N h\rho_j k_0(x_i,x_j)U_m(x_j) ,\\
Z_m(x_i) &= \sum _{j=0}^N h\rho_j k_1(x_i,x_j)Y_ m(x_j) ,\;  i=0,...,N.
\end{split}
\end{equation}


In order to get the error estimates for the approximate solution of the problem \eqref{eq3} and its derivatives on the grid we need some following auxiliary results.  
\begin{Proposition}\label{prop1}
Assume that the function $f(t,u,y,v,z)$ has all continuous partial derivatives up to second order in the domain $\mathcal{D}_M$, the kernel functions $k_0(x,t), k_1(x,t)$ have all continuous partial derivatives up to second order in the  square $[0,1] \times [0,1]$ and the function $\varphi (t)$ has continuous derivatives up to second order. Then for the functions $\psi_m(x), u_m(x), y_m(x),  v_m(x), z_m(x), m=0,1,...,$ constructed by the iterative method \eqref{iter1c}-\eqref{iter3c} we have  
$\psi_m(x) \in  C^2 [0, 1]$, $u_m(x) \in C^6 [0, 1]$, $ y_m(x), v_m(x), z_m(x) \in C^2 [0, 1]$.
\end{Proposition}
This proposition is obvious.
%

\begin{Proposition}\label{prop2}

For any function $\psi (x) \in C^2[0, 1]$ there holds the estimate
\begin{equation}\label{eq:prop2}
\begin{split}
&\int_0^1 G (x_i,t) \psi (t) dt = \sum _{j=0}^N h\rho_j G(x_i,t_j)\psi(t_j) +O(h^2),\\
& \int_0^1 G (\varphi (x_i),t) \psi (t) dt = \sum _{j=0}^N h\rho_j G(\varphi (x_i),t_j)\psi(t_j) +O(h^2).
\end{split}
\end{equation}
\end{Proposition}
\noindent {\bf Proof.}
 The above estimates are obvious in view of the error estimate of the compound trapezium formula.
 

\begin{Proposition}\label{prop3}
Under the assumptions of Proposition \ref{prop1},  for any $m=0,1,...$ there hold the estimates
\begin{equation}\label{eq:prop3a}
\|\Psi_m -\psi_m  \|= O(h^2),\; \|U_m -u_m  \|=O(h^2),\; \|Y_m -y_m  \|=O(h^2),
\end{equation}
\begin{equation}\label{eq:prop3b}
\begin{split}
  \|V_m-v_m  \|=O(h^2), \; \|Z_m -z_m  \|=O(h^2).
\end{split}
\end{equation}
 where $ \|.\|=  \|.\|_{\bar{\omega}_h}$ is the max-norm of function on the grid $\bar{\omega}_h$.
\end{Proposition}
\noindent {\bf Proof.}
We prove the proposition by induction. For $m=0$ we have immediately $\|\Psi_0 -\psi_0  \|= 0$. Next, by the first equation in \eqref{iter2c} and Proposition \ref{prop2} we have
\begin{equation}
u_0(x_i)=\int_0^1 G (x_i,t) \psi_0 (t) dt = \sum _{j=0}^N h\rho_j G(x_i,t_j)\psi_0(t_j)+O(h^2)
\end{equation}
for any $i=0,...,N$ .
On the other hand, in view of  the first equation in \eqref{iter2d}  we have
\begin{equation}
U_0(x_i)= \sum _{j=0}^N h\rho_j G_0(x_i,t_j)\Psi_0(t_j).
\end{equation}
Therefore, $|U_0(t_i)- u_0(t_i)|= O(h^2)$ because $\Psi_0(t_j)=\psi_0(t_j)=f(t_j,0,0,0,0)$. Consequently, $\|U_0 -u_0  \|=O(h^2) $.\\
Similarly, we have 
\begin{equation}
\|Y_0 -y_0  \|=O(h^2).
\end{equation}
Next, by the trapezium formula we have
\begin{align*}
v_0(x_i)= \int_0^1 k_0(x_i,t)u_0(t) dt =  \sum _{j=0}^N h\rho_j k_0(x_i,t_j)u_0(t_j) +O(h^2),
\end{align*}
while by the third equation in \eqref{iter2d} we have
\begin{align*}
V_0(x_i) = \sum _{j=0}^N h\rho_j k_0(x_i,t_j)U_0(t_j) ,\;  i=0,...,N.
\end{align*}
Therefore, 
\begin{align*} 
\Big | V_0(x_i) - v_0(x_i)   \Big | &= \Big | \sum _{j=0} ^{N} h\rho_j k_0(x_i,t_j) (U_0(t_j) -u_0(t_j) )  \Big | +O(h^2) \\
&\le  \sum _{j=0} ^{N} h\rho_j |k_0(x_i,t_j)| |U_0(t_j) -u_0(t_j) |   +O(h^2) \\
&\le C h^2  \sum _{j=0} ^{N} h\rho_j |k_0(x_i,t_j)| + O(h^2) \\
&\le CC_1 h^2 \sum _{j=0} ^{N} h\rho_j + O(h^2) =  O(h^2)
\end{align*} 
because $|U_0(t_j) -u_0(t_j) | \le Ch^2, \; |k_0(x_i,t_j)| \le C_1$, where $C, C_1$ are some constants. Thus, $\|V_0-v_0 \| =O(h^2)$.
Analogously, we have $\|Z_0-z_0\| =O(h^2) $.

Now suppose that \eqref{eq:prop3a} and \eqref{eq:prop3b} are valid for $m \ge 0$. We shall show that these estimates are valid for $m+1$.
By the Lipschitz condition of the function $f$ and the estimates \eqref{eq:prop3a} and \eqref{eq:prop3b} it is easy to obtain the estimate $$\|\Psi_{m+1} -\psi_{m+1}  \|= O(h^2).$$ 
Now from the first equation in \eqref{iter2c} by Proposition \ref{prop2} we have
\begin{equation*}
u_{m+1}(x_i)=\int_0^1 G (x_i,t) \psi_{m+1} (t) dt = \sum _{j=0}^N h\rho_j G(x_i,x_j)\psi_{m+1}(x_j)+O(h^2).
\end{equation*}
On the other hand by the first formula in \eqref{iter2d} we have
\begin{equation*}
U_{m+1}(x_i) = \sum _{j=0}^N h\rho_j G(x_i,x_j)\Psi_{m+1}(x_j).
\end{equation*}
From this equality and the above estimates  we obtain the estimate
$$\|U_{m+1} -u_{m+1}  \|=O(h^2).
$$
Similarly, we obtain
$$\|Y_{m+1} -y_{m+1}\|=O(h^2),\; \|V_{m+1} -v_{m+1}\|=O(h^2),\; \|Z_{m+1} -z_{k+1}\|=O(h^2).
$$
Thus, by induction we have proved the proposition.
Now combining Proposition \ref{prop3} and Theorem \ref{thm3} results in the following theorem.
\begin{Theorem}\label{thm4}
 Assume that all the conditions of Theorem \ref{thm1} and Proposition \ref{prop1} are satisfied. Then, for the approximate solution of the problem \eqref{eq3} obtained by the discrete iterative method on the uniform grid with grid size $h$ there hold the estimates for the total error
\begin{equation}\label{eqthm4}
\|U_m-u\| \leq M_0 p_md +O(h^2). 
\end{equation}
\end{Theorem}
\noindent {\bf Proof.}
The first above estimate is easily obtained if representing
\begin{equation*}
U_m(t_i)-u(t_i)= (u_m(t_i)-u(t_i))+(U_m(t_i)-u_m(t_i))
\end{equation*}
and using  the  estimate in Theorem \ref{thm3} and the second estimate in \eqref{eq:prop3a}. Thus, the theorem is proved.

\section{Examples}

\noindent {\bf Example 1.} Consider the problem 
\begin{equation}\label{exam1}
\begin{split}
u^{(4)}(x)& = \pi^4 \sin (\pi x)-\frac{1}{2}\sin^2 (\pi x)-\frac{1}{2}\sin^2 (\frac{\pi}{2}x)+\frac{1}{2}u^2(x)+\frac{1}{2}u^2(\frac{x}{2})\\
&-\frac{8}{3 \pi}\int_0^1 e^x  \sin(\pi t) u(t)dt + \int_0^1 e^x  \sin(\pi t) u(\frac{t}{2})dt,\\
u(0)&=0 , \; u(1)=0, \; u''(0)=0 , \; u''(1)=0.
\end{split}
\end{equation}
In this problem
\begin{align*}
&k_0(x,t)=k_1(x,t)= e^x \sin (\pi t),\; \varphi (t)=\frac{t}{2},   \\
&f(x,u,y,v,z)= \pi^4 \sin (\pi x)-\frac{1}{2 }\sin^2 (\pi x)-\frac{1}{2}\sin^2 (\frac{\pi}{2}x)+\frac{1}{2}u^2+\frac{1}{2}y^2-\frac{8}{3 \pi}v + z.
\end{align*}
So,
 \begin{equation*}
K_0=K_1= \max_{0\leq x\leq 1} \int_0^1 |k_0(x,s)|ds =\frac{2 e}{\pi}.
\end{equation*}
It is easy to verify that the function $u= \sin (\pi x)$ is the exact solution of the problem.  In the domain $\mathcal{D}_M$  defined by 
\begin{align*} 
\mathcal{D}_M = \{(x,u,y,v,z)\; | \; &0\leq x\leq 1,\; |u|\leq M_0 M, \;|y| \le M_0M,\;\\
& |v| \le K_0M_0M,\; |z|\leq K_1M_0 M \}
\end{align*}
we have
\begin{align*}
|f(x,u,y,v,z)| \le  \pi ^4 + 1+ M_0^2M^2 +(\frac{8}{3 \pi}+1)\frac{2e}{\pi}M_0M.
\end{align*}
Taking into account that $M_0=\frac{5}{384}$ it is possible to verify that for $M=105$ all the conditions of Theorem \ref{thm1} are satisfied with
$ L_0 = 1.3672, L_1= 1.4714, L_2=0.8488, L_3=1, q=0.0773$.
Therefore, the problem has a unique solution $u(x)$ satisfying the estimates $|u(x)| \le 1.3672 $. These theoretical estimates are somewhat greater than the exact estimates $|u(x)| \le 1 $.\par

Below we report the numerical results by the discrete iterative method  \eqref{iter1d}-\eqref{iter3d}  for the problem. In Tables \ref{table:1} and \ref{table:2} we use the notation $Error =\|U_m-u \|$, where $u$ is the exact solution of the problem.
\begin{table}[ht!]
\centering
\caption{The convergence in Example 1 for the stopping criterion $\|U_m-u \|\le h^2$}
\label{table:1}
\begin{tabular}{cccc}
\hline 
$N$ &	$h^2$ &	$m$	&$Error$ \\
\hline 
50	& 4.0000e-04 &	2 &	1.1564e-04\\
100	& 1.0000e-04 &	3&	2.2752e-06\\
150	&4.4444e-05	 & 3 & 1.9519e-06\\
200	&2.5000e-05	& 3	& 1.8386e-06\\
300	&1.1111e-05	& 3	& 1.7575e-06\\
400& 6.2500e-06 & 3	& 1.7292e-06\\
500	&4.0000e-06	&3	& 1.7160e-06\\
800	&1.5625e-06&	4&	3.4384e-08\\
1000&	1.0000e-06	&4	& 3.1098e-08\\
\hline 
\end{tabular} 
\end{table}
 It is interesting to notice that if taking the stopping criterion $\|\Psi_m-\Psi_{m-1} \|\le 10^{-9}$  instead of $\|U_m-u \|\le h^2$ then we obtain better accuracy of the approximate solution with more iterations. See Table \ref{table:2}. 

\begin{table}[ht!]
\centering
\caption{The convergence in Example 1 for the stopping criterion $\|\Psi_m-\Psi_{m-1} \|\le 10^{-9}$}
\label{table:2}
\begin{tabular}{cccc}
\hline 
$N$ &	$h^2$ &	$m$	&$Error$ \\
\hline 

50	& 4.0000e-04 &	6 &	2.3139e-06\\
100	& 1.0000e-04 &	6&	5.8292e-07\\
150	&4.4444e-05	& 6	&2.5941e-07\\
200	&2.5000e-05	& 6	&1.4600e-07\\
300	&1.1111e-05	&6	&6.4911e-08\\
400&	6.2500e-06&	6 &3.6519e-08\\
500	&4.0000e-06	&6	& 2.3376e-08\\
800	&1.5625e-06&	6&	9.1351e-09\\
1000&	1.0000e-06	&6	&5.8485e-09\\
\hline 
\end{tabular} 
\end{table}


From Table \ref{table:2} we see that the accuracy of the approximate solution is much better than $O(h^2)$.\\

\noindent {\bf Example 2.}  Consider the nonlinear fourth order BVP
\begin{equation}\label{exam2}
\begin{split}
u^{(4)}(x)& =1+\pi^2 \sin (\pi x)+2u^2(x)+2 u^2(\frac{x}{2})+e^{-u^2(x)}\\
&+3 \Big (\int_0^1 e^x  \sin(\pi t ) u(t)dt\Big )^2 \Big( \int_0^1 e^x  \sin(\frac{\pi tx}{2}) u(\frac{t}{2})dt\Big)^2,\\
u(0)&=0 , \; u(1)=0, \; u''(0)=0 , \; u''(1)=0.
\end{split}
\end{equation}
In this problem
\begin{align*}
&k_0(x,t)= e^x \sin (\pi t),\; k_1(x,t)= e^x \sin(\frac{\pi tx}{2}),\;  \varphi (t)=\frac{t}{2},   \\
&f(x,u,y,v,z)= 1+\pi^2 \sin (\pi x)+2u^2+2y^2+e^{-u^2}+3v^2 z^2
\end{align*}
It is possible to verify that $K_1 =\max_{0\leq x\leq 1} \int_0^1 |k_1(x,t)|dt = \frac{2e}{\pi}$
and for $M=12.5$ all the assumptions of Theorem \ref{thm2} are satisfied.
Therefore, the problem \eqref{exam2} has a unique positive solution satisfying $0 \le u(x) \le  0.1628$. \\
Using the discrete iterative method \eqref{iter1d}-\eqref{iter3d} on the grid with grid step $h=0.01$ and the stopping criterion $\|\Psi_m-\Psi_{m-1} \|\le 10^{-9}$ we found an approximate solution after $5$ iterations. The graph of this approximate solution is depicted in Figure 1.\\

\begin{figure}[ht]
\begin{center}
\includegraphics[height=6cm,width=9cm]{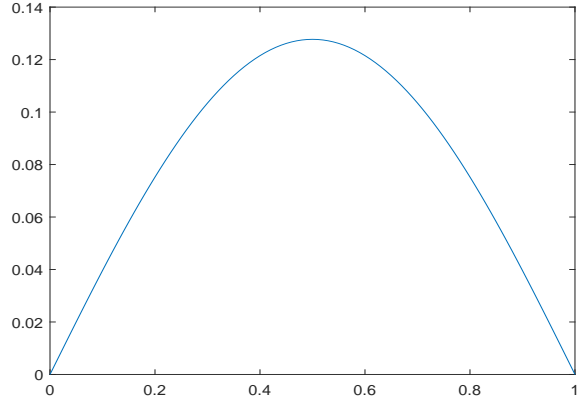}
\caption{The graph of the approximate solution in Example $2$. }
\label{fig1}
\end{center}
\end{figure}

\noindent {\bf Example 3.}  Consider  the nonlinear fourth order BVP with non-homogeneous boundary conditions
\begin{equation}\label{exam3}
\begin{split}
u^{(4)}(x)& =24-x^4(1-x)^4-(\frac{x}{3})^4 (1- \frac{x}{3})^4 - \frac{e^{x} }{60}\\
& - \frac{4(\pi^4 - 3\pi^2 + 12)}{81\pi^5} x^2 (1-x)^2 \sin (\pi x) + u^2(x)  \\
&+u^2(\frac{x}{3})+\int_0^1 k_0(x,t) u(t) dt +u(x) \int_0^1 k_1(x,t) u(\frac{t}{3})dt,\\
u(0)&=0 , \; u(1)=0, \; u''(0)=2 , \; u''(1)=2,
\end{split}
\end{equation}
where
\begin{equation*}
k_0(x,t)=e^{x}t,\; k_1(x,t)=\sin (\pi x)\sin(\pi t).
\end{equation*}
This problem has the exact solution $u(x)=x^2 (1-x)^2$.\\
Using the discrete iterative method \eqref{iter1d}, \eqref{iter2da},\eqref{iter3d} we obtain the results of convergence reported in Table \ref{table3}.
\begin{table}[ht!]
\centering
\caption{The convergence in Example 3 for the stopping criterion $\|\Psi_m-\Psi_{m-1} \|\le 10^{-9}$}
\label{table3}
\begin{tabular}{cccc}
\hline 
$N$ &	$h^2$ &	$m$	&$Error$ \\
\hline 

50	& 4.0000e-04 &	5 &	1.0091e-04\\
100	& 1.0000e-04 &	5&	5.2227e-05\\
150	&4.4444e-05	& 5	& 1.1212e-05 \\
200	&2.5000e-05	& 5	& 6.3068e-06 \\
300	&1.1111e-05	& 5	& 2.8030e-06 \\
400&	6.2500e-06&	5 & 1.5767e-06 \\
500	&4.0000e-06	& 5	& 1.0091e-06 \\
800	&1.5625e-06& 5 & 3.9417e-07 \\
1000&	1.0000e-06	&6	&2.5227e-07 \\
\hline 
\end{tabular} 
\end{table}

\noindent {\bf Remark 2.} The iterative method \eqref{iter1d}- \eqref{iter3d} was proved to be convergent under the conditions of Theorem \ref{thm1}. When these conditions are not satisfied Theorem \ref{thm3}
 does not guarantee the convergence of the iterative method. But in many cases the iterative method remains convergent. Below we show an example, where although the right-hand function has a weak singularity at $x=0$ the method  converges. \\
 
\begin{figure}[ht]
\begin{center}
\includegraphics[height=6cm,width=9cm]{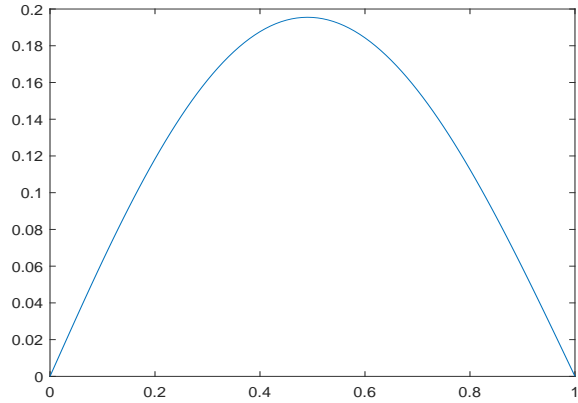}
\caption{The graph of the approximate solution in Example $2$. }
\label{fig2}
\end{center}
\end{figure}

\noindent {\bf Example 4.}
Consider Example 2 with the new right-hand function of the form
 \begin{equation*}
f(x,u,y,v,z)= \frac{1}{\sqrt x}(1+\pi^2 \sin (\pi x)+2u^2+2y^2+e^{-u^2}+3v^2 z^2).
 \end{equation*}
Notice that in this case, since the function $\psi (t)$ has a weak singularity at $x=0$,  the integral
$\int _0^1 G(x,s)\psi(s) ds $ exists. Consequently, the continuous iterative method \eqref{iter1c}- \eqref{iter3c} can be carried out. In the discrete iterative method \eqref{iter1d}- \eqref{iter3d}, to avoid computing $\Psi_m(x_1)$ for $x_1=0$ we set $\Psi_m(x_1)=0$ due to $G(x_i,x_1)=0$. The numerical experiments show that the discrete iterative method converges, and with $TOL=10^{-9}$ the iterative process stops after $6$ iterations. The graph of the approximate solution is depicted in Figure \ref{fig2}.

\section{Conclusion}
In this paper we have established the existence, uniqueness and positivity of the solution for a fourth order nonlinear functional integro-differential equation with the Navier boundary conditions and proposed an iterative method at both continuous and discrete levels for finding the solution. The second order of accuracy of the discrete method has been proved. Some examples, where the exact solution is known and is not known, demonstrate the validity of the obtained theoretical results and the efficiency of the iterative method. \par
The method used in this paper with appropriate modifications can be applied to nonlinear functional integro-differential equations of any order with other boundary conditions and more complicated nonlinear terms. This is the direction of our research in the future. 


\end{document}